\documentclass[12pt,leqno]{article}
\usepackage{latexsym,amsfonts,epsfig}
\usepackage{amsmath,amssymb,amsthm}
\usepackage[notref, notcite]{}
\usepackage[colorlinks,linkcolor=red,citecolor=blue,urlcolor=blue]{hyperref}
\usepackage[margin=0.5in]{geometry}
\usepackage{graphicx}
\input epsf
\newtheorem{theorem}{Theorem}[section]
\newtheorem{lemma}{Lemma}[section]

\RequirePackage[colorlinks,citecolor=blue,urlcolor=blue]{hyperref}

\begin{document}
\title{Inference via the Skewness-Kurtosis Set}
\author{ Chris A.J. Klaassen and Bert van Es \\
Korteweg-de Vries Insitute for Mathematics \\
University of Amsterdam}

\maketitle



\begin{abstract}
Kurtosis minus squared skewness is bounded from below by 1, but for unimodal distributions this parameter is bounded by 189/125.
In some applications it is natural to compare distributions by comparing their kurtosis-minus-squared-skewness parameters.
The asymptotic behavior of the empirical version of this parameter is studied here for i.i.d. random variables.
The result may be used to test the hypothesis of unimodality against the alternative that the kurtosis-minus-squared-skewness parameter is less than 189/125.
However, such a test has to be applied with care, since this parameter can take arbitrarily large values, also for multimodal distributions.
Numerical results are presented and for three classes of distributions the skewness-kurtosis sets are described in detail.
\end{abstract}

\section{Introduction\label{intro}}

Skewness and kurtosis play an important role in applied research; see e.g. \cite{Cristelli} (complex dynamics), \cite{BusseJelly} (turbulence), \cite{Martin} (health care), \cite{Dunnwald} (medicine, deep learning), \cite{Hughes} (ocean modelling), \cite{Karagiorgis} (cryptocurrencies).
Comparing their data to the normal distribution these researchers study the standardized third and fourth moments, which means that they study the skewness and kurtosis.
If $X$ is a random variable with finite fourth moment, mean $\mu$ and variance $\sigma^2$, then the skewness $\tau$, kurtosis $\kappa$ and excess kurtosis $\kappa_{\rm e}$ are defined as
\begin{equation*}
\tau = \frac{ E \left( (X - \mu)^3 \right)} {\sigma^3},\quad \kappa = \frac{ E \left(  (X - \mu)^4 \right)} {\sigma^4}, \quad \kappa_{\rm e} = \kappa -3.
\end{equation*}
Since normal distributions have kurtosis 3, the excess kurtosis is a measure for the deviation from Gaussianity.
In practice one often calls the excess kurtosis simply kurtosis.
Consider the plane with skewness on the horizontal ($x$-) and kurtosis on the vertical ($y$-)axis.
Pearson \cite{Pearson} published the inequality
\begin{equation}\label{Pearson}
\kappa \geq \tau^2 + 1,
\end{equation}
which holds for all distributions; see (7)--(9) of \cite{Mokveld}.
This means that all points $(\tau,\kappa)$ are on or above the parabola $\{ (x,y) \in {\mathbb R}^2 \, | \, y = x^2 +1 \}$.
In Subsection \ref{sk} we prove that all points on or above this parabola can be attained by choosing the underlying distribution appropriately.

For unimodal distributions the Pearson inequality was sharpened to
\begin{equation}\label{Philip}
\kappa \geq \tau^2 + 189/125
\end{equation}
by \cite{Mokveld}; a distribution is unimodal if its  distribution function is convex--concave.
In Subsection \ref{skunimodal} we show that for unimodal distributions all points $(\tau,\kappa)$ are above the parabola
$\{ (x,y) \in {\mathbb R}^2 \, | \, y = x^2 +189/125 \}$, except for two points that are at this parabola.
Actually, we sharpen and optimize inequality (1.2) there by replacing 189/125 by a complicated function of $\tau$.

For symmetric unimodal distributions the inequality
\begin{equation}\label{Bert}
\kappa \geq \tau^2 + 9/5
\end{equation}
holds; see Remark 2.1 and Table 1 of \cite{Mokveld}.
In Subsection \ref{sksunimodal} we show that for symmetric unimodal distributions the skewness-kurtosis set of all points $(\tau,\kappa)$ equals the half line $\{ (x,y) \in {\mathbb R}^2 \, | \, x=0,\, y \geq 9/5 \}$.

In \cite{Karagiorgis} one starts with daily returns on cryptocurrencies.
On a weekly basis, sample skewness and sample kurtosis of these daily returns are computed.
In this way many skewness-kurtosis points in the plane are obtained.
By a regression approach it is argued in \cite{Karagiorgis} that the regression line $\kappa = \tau^2 + 189/125$ (see (\ref{Philip})) fits these data points better than the regression lines $\kappa = \tau^2 + 1$ corresponding to (\ref{Pearson}) and $\kappa = \tau^{4/3}$ suggested by \cite{Cristelli}.
Another approach would have been to test the hypothesis $\Delta = \kappa - \tau^2 \geq 189/125$.
Our main result is Theorem \ref{general}, which presents the asymptotic behavior of the empirical version of $\Delta$ under distributions with finite eighth moments.

Inequality (\ref{Philip}) implies that distributions with $\Delta < 189/125$ are necessarily multimodal.
However, there are multimodal distributions with arbitrarily large positive values of $\Delta$; see (\ref{Deltalarge}).

In Section \ref{estimator} we introduce the test statistic based on the empirical skewness and kurtosis.
Since it turns out that the boundary inflated uniform distribution plays a crucial role, we consider the test statistic under this distribution in more detail in Section \ref{boundarinflateduniform}.
It is also observed that for several multimodal distributions the bound 189/125 is surpassed.
In Section \ref{multimodality} this is investigated for certain mixtures of two normal distributions and of two exponential distributions.
For three classes of distributions, namely the class of all distributions, of all unimodal distributions and of all symmetric unimodal distributions, the skewness-kurtosis set is given a complete description in Section \ref{sksets}.
Finally, Section \ref{proofs} contains proofs of Theorems \ref{BIU} and \ref{MN}.

\section{An empirical moments estimator\label{estimator}}

Let $X, X_1, \dots, X_n$ be i.i.d. random variables with mean $\mu$, variance $\sigma^2 >0$, skewness $\tau = \sigma^{-3} E((X-\mu)^3)$,
kurtosis $\kappa= \sigma^{-4} E((X-\mu)^4)$ and finite eighth moment.
Let
\begin{equation*}
T_n = \frac{ \tfrac 1n \sum_{i=1}^n \left(X_i - {\bar X}_n \right)^4 }
{\left[ \tfrac 1n \sum_{i=1}^n \left(X_i - {\bar X}_n \right)^2 \right]^2} \,
- \, \frac{\left[ \tfrac 1n \sum_{i=1}^n \left(X_i - {\bar X}_n \right)^3 \right]^2}
{\left[ \tfrac 1n \sum_{i=1}^n \left(X_i - {\bar X}_n \right)^2 \right]^3}, \quad {\bar X}_n = \tfrac 1n \sum_{i=1}^n X_i,
\end{equation*}
be the empirical version of $\Delta = \kappa - \tau^2$, which estimates $\Delta$.
Since both $\Delta$ and $T_n$ are location and scale invariant, the behavior of $T_n - \Delta$ is the same for all members of the location-scale family of $X$.
Therefore we assume $\mu=0$ and $\sigma^2=1$, without loss of generality.

We denote the $k$th (reduced) moment of $X$ by
\begin{equation*}
\nu_k = E\left( \left( \frac{X-\mu}{\sigma} \right)^k \right) = E\left( X^k \right), \quad k \leq 8.
\end{equation*}
By the central limit theorem and the boundedness of the eighth reduced moment $\nu_8$ of $X$
\begin{equation*}
S_{k,n} = \tfrac 1n \sum_{i=1}^n \left[ X_i^k - \nu_k \right] = {\cal O}_p\left( \frac 1{\sqrt n} \right), \quad k \leq 4,
\end{equation*}
holds.
Hence we have
\begin{eqnarray*}
\lefteqn{ \tfrac 1n \sum_{i=1}^n \left(X_i - {\bar X}_n \right)^k
= \tfrac 1n \sum_{i=1}^n \left( X_i^k -k {\bar X}_n X_i^{k-1} \right) + {\cal O}_p \left(\tfrac 1n \right) } \\
&& = \nu_k + S_{k,n} -k \nu_{k-1} {\bar X}_n + {\cal O}_p \left(\tfrac 1n \right), \quad k \leq 4,
\end{eqnarray*}
and
\begin{eqnarray*}
\lefteqn{ T_n - \Delta = \left[1+ {\cal O}_p\left( \tfrac 1{\sqrt n} \right) \right]^{-3} } \\
&& \times \left\{ \left[1+ S_{2,n} + {\cal O}_p \left(\tfrac 1n \right) \right]
\left[\nu_4 + S_{4,n} - 4 \nu_3 {\bar X}_n + {\cal O}_p \left(\tfrac 1n \right) \right] \right. \\
&& \quad \left. - \left[\nu_3 + S_{3,n} -3{\bar X}_n + {\cal O}_p \left(\tfrac 1n \right) \right]^2
     - \left( \kappa - \tau^2 \right) \left[1+ S_{2,n} + {\cal O}_p \left(\tfrac 1n \right) \right]^3 \right\} \nonumber \\
&& = \frac 1{\left[1+ {\cal O}_p\left( \frac 1{\sqrt n} \right) \right]^3} \left\{ \nu_4 + S_{4,n} - 4 \nu_3 S_{1,n} + \nu_4 S_{2,n} - \nu_3^2 - 2\nu_3 S_{3,n} \right. \\
&& \quad \quad \left. + 6 \nu_3 S_{1,n} - \left( \nu_4 - \nu_3^2 \right) \left[1 + 3 S_{2,n} \right] + {\cal O}_p \left(\tfrac 1n \right) \right\} \\
&& = \frac { 2 \nu_3 S_{1,n} - \left( 2 \nu_4 - 3 \nu_3^2 \right) S_{2,n} - 2\nu_3 S_{3,n} + S_{4,n} + {\cal O}_p \left(\tfrac 1n \right) }
{ \left[1+ {\cal O}_p\left( \frac 1{\sqrt n} \right) \right]^3 } \nonumber \\
&& = 2 \nu_3 S_{1,n} - \left( 2\nu_4 - 3 \nu_3^2 \right) S_{2,n} - 2\nu_3 S_{3,n} + S_{4,n} + {\cal O}_p \left(\tfrac 1n \right). \nonumber
\end{eqnarray*}
We obtain
\begin{eqnarray*}
\lefteqn{ \sqrt{n} \left( T_n - \Delta \right) } \\
&& = \frac 1{\sqrt n} \sum_{i=1}^n \left[ 2 \nu_3 X_i - \left( 2 \Delta - \nu_3^2 \right)( X_i^2 - 1)
- 2 \nu_3 ( X_i^3 - \nu_3) + ( X_i^4 - \nu_4 ) \right] \\
&& \hspace{10em} + {\cal O}_p \left(\tfrac 1{\sqrt n} \right).
\end{eqnarray*}
Straightforward computation shows that this implies

\begin{theorem}\label{general}
Let $X, X_1, \dots, X_n$ be i.i.d. random variables with mean $\mu$, variance $\sigma^2 >0$, skewness $\tau = \sigma^{-3} E((X-\mu)^3)$,
kurtosis $\kappa = \sigma^{-4} E((X-\mu)^4)$ and finite eighth moment, and let $\nu_k = E\left( \left( \frac{X-\mu}{\sigma} \right)^k \right)$ be the $k$-th reduced moment.
With $\Delta = \kappa - \tau^2 = \nu_4 - \nu_3^2$
\begin{equation*}
\sqrt{n} \left( T_n - \Delta \right) \stackrel{\cal D} {\longrightarrow} {\cal N}(0, \varsigma^2), \quad {\rm as~} n \to \infty
\end{equation*}
holds, where the asymptotic variance $\varsigma^2$ is given by
\begin{eqnarray}\label{Aemeh}
\lefteqn{ \varsigma^2 = 4 \Delta^3 - \Delta^2 - \left( 3 \nu_3^4 + 16 \nu_3^2 - 8 \nu_3 \nu_5 + 4 \nu_6 \right) \Delta } \\
&& + \nu_3^6- 4 \nu_3^4 - 4 \nu_3^3 \nu_5 + 2 \nu_3^2 \left( 2 + 3 \nu_6 \right) + 4 \nu_3 \left( \nu_5 - \nu_7 \right) + \nu_8. \nonumber
\end{eqnarray}
\end{theorem}

This result can be obtained also by applying Example 2 or 3 of \cite{samplemoments} and the multivariate $\delta$-technique.
Actually, our proof follows similar lines, but uses the peculiarities in the structure of $T_n$.

To test the null hypothesis of unimodality, which implies $\Delta \geq 189/125$ or equivalently $\tau^2 - \kappa_{\rm e}\leq 186/125$, we consider the random variable $X$ that has point mass 1/2 at $-\sqrt{3/5}$ and with probability 1/2 behaves like a uniformly distributed random variable on $(-\sqrt{3/5}, 3 \sqrt{3/5})$.
$X$ has a boundary inflated uniform distribution with point mass 1/2, mean 0 and variance 1.
By \cite{Mokveld} it was shown that for unimodal distributions equality in $\tau^2 - \kappa_{\rm e}\leq 186/125$ holds for unimodal distributions if and only if the underlying distribution is boundary inflated uniform with point mass 1/2.

In order to apply Theorem \ref{general} for such boundary inflated uniform distributions we need to compute (recall that $T_n$ and $\Delta$ are location and scale invariant)
\begin{eqnarray*}
\lefteqn{ \nu_k = E \left( X^k \right) = \frac 12 \left( - \sqrt{\frac 35} \right)^k +
\frac 12 \sqrt{\frac 53} \ \frac 1{4(k+1)} \left[ \left(3 \sqrt{\frac 35} \right)^{k+1} - \left( -\sqrt{\frac 35} \right)^{k+1} \right] } \\
&& \hspace{12em} = \left( \sqrt{\frac 35} \right)^k  \frac {3^{k+1} +(4k+5)(-1)^k }{8(k+1)}, \quad k \leq 8,
\end{eqnarray*}
which means
\begin{eqnarray*}
\lefteqn{ \nu_1 =0, \quad \nu_2 = 1, \quad \nu_3 = \frac 65 \sqrt{\frac 35}, \quad \nu_4 = \frac {297}{125}, \quad
\nu_5 = \frac{132}{25}\sqrt{\frac 35} = \frac {22}5 \nu_3, \nonumber } \\
&& \hspace {6em} \nu_6 = \frac{7479}{875}, \quad \nu_7 = \frac{2754}{125}\sqrt{\frac 35} = \frac{459}{25} \nu_3, \quad
\nu_8 = \frac {4437}{125}.
\end{eqnarray*}
Straightforward but rather tedious computation shows that these values yield the limit variance
\begin{equation}\label{Aemej}
\varsigma^2 = \frac { 2^{16} \times 3^2 } { 5^7 \times 7 } \approx 1.0785353143, \quad \varsigma \approx 1.0385255482.
\end{equation}
Together with Theorem 1 from \cite{Mokveld}, Theorem \ref{general} and (\ref{Aemej}) imply that for i.i.d. random variables $X_1, \dots, X_n$ with a distribution with finite eighth moment, one may test the null hypothesis of unimodality (and hence $\Delta \geq 189/125$) versus the alternative hypothesis of $\Delta < 189/125$ (and hence multimodality) at significance level $\alpha$ by rejecting whenever for large values of $n$
\begin{equation}\label{test}
T_n - \frac{189}{125} \leq - \frac {\varsigma}{\sqrt n} \ \Phi^{-1}(1-\alpha), \quad \varsigma \approx 1.0385255482,
\end{equation}
holds.
Although all unimodal distributions satisfy $\Delta \geq 189/125$ (see \cite{Mokveld}, but also Subsection 5.2), there are multimodal distributions as well satisfying this inequality, as one may see from the following simple argument.

Let $U$ and $V$ be independent random variables, $U$ with a uniform distribution on the unit interval and $V$ with a Bernoulli distribution taking the values 2 and 3 each with probability 1/2.
Let $W$ be equal to $V$ with probability $p$ and equal to $U$ with probability $1-p$.
Denoting the skewness of $W$ by $\tau(p)$ and its kurtosis by $\kappa(p)$ we observe that $\Delta(p)= \kappa(p) - \tau^2(p)$ is a continuous function of $p \in [0,1]$.
By e.g. Table 1 from \cite{Mokveld} we have $\Delta(0) = 9/5 > 189/125 > 1 = \Delta(1)$.
Finally, note that $W$ has a multimodal distribution for $0<p<1$.

\section{Boundary inflated uniform, BIU($p$)} \label{boundarinflateduniform}

The location scale family of boundary inflated uniform distributions may be represented by the distribution that has probability $1-p$ for zero and conditionally on being positive it renders a uniform distribution on $(0,1]$.
The normalized version of this distribution, i.e. with mean 0 and variance 1, is the distribution of the random variable
\begin{equation*}
X = \left\{ \begin{array}{lcr}
            - \sqrt{ \frac{3p}{4-3p}}    &                        & 1-p \\
                                         & {\rm with~probability} &     \\
            {\rm uniform~on} \left[ - \sqrt{ \frac{3p}{4-3p}}, \frac {2-p}p \sqrt{ \frac{3p}{4-3p}} \right] & & p.
            \end{array}
    \right.
\end{equation*}
It can be represented as
\begin{equation*}
X = \frac 2p \sqrt{ \frac{3p}{4-3p}} \left[BU - \frac p2 \right]
\end{equation*}
with $B$ and $U$ independent, $B$ Bernoulli distributed with success probability $p$ and $U$ uniform on the unit interval.

The following theorem gives the values of the (reduced) moments, skewness, kurtosis, the test function $\Delta$ and the limit variance $\varsigma^2$ of the test statistic for the distribution of $X$.
\begin{theorem}\label{BIU}
For the {\rm BIU($p$)} distribution the moments are
\begin{eqnarray*}
\lefteqn { \nu_k = E X^k = \left[ \sqrt{ \frac{3p}{4-3p}} \right]^k \left\{ (-1)^k (1-p) +
                    \frac{(2-p)^{k+1}-(-p)^{k+1}}{2(k+1)p^{k-1}} \right\} \nonumber } \\
&& \nu_1 = 0,\  \nu_2 = 1,\ \nu_3 = \frac{6(1-p)^2}{p(4-3p)} \sqrt{ \frac{3p}{4-3p}},\
                            \nu_4 = \frac {9(16 - 40p + 40p^2 -15 p^3)}{5p(4-3p)^2}  \nonumber \\
&& \nu_5 = \frac {12(1-p)^2(4 - 4p + 3p^2)}{p^2(4-3p)^2} \sqrt{ \frac{3p}{4-3p}}, \\
&& \nu_6 = \frac{27(64 - 224p + 336p^2 - 280 p^3 + 140 p^4 - 35p^5)}{7p^2(4-3p)^3}, \nonumber \\
&& \nu_7 = \frac{54(1-p)^2(8-16p+16p^2-8p^3+3p^4)}{p^3(4-3p)^3} \sqrt{ \frac{3p}{4-3p}}, \nonumber \\
&& \nu_8 = \frac {9(256 - 1152p + 2304p^2 - 2688p^3 + 2016p^4 - 1008p^5 + 336p^6 -63p^7)}{p^3(4-3p)^4} \nonumber
\end{eqnarray*}
the skewness is given by
\begin{equation*}
\tau(p) = \nu_3 = EX^3 = \frac{6(1-p)^2}{p(4-3p)} \sqrt{ \frac{3p}{4-3p}},
\end{equation*}
the kurtosis by
\begin{equation*}
\kappa(p) = \nu_4 = EX^4 = \frac {9(16 - 40p + 40p^2 -15 p^3)}{5p(4-3p)^2}
\end{equation*}
and the test function, kurtosis minus squared skewness, by
\begin{equation}\label{testBIU}
\Delta(p) = \kappa(p) - \tau^2(p) = \frac {9(4 + 32p - 80 p^2 + 60 p^3 - 15 p^4)}{5p(4-3p)^3}.
\end{equation}
Finally the limit variance in (\ref{Aemeh}) of the test statistic equals
\begin{equation}\label{statBIU}
\varsigma^2(p) = \frac{1152(2446-11229p+21744p^2-22080p^3+11640p^4-2520p^5)}{875p^3(4-3p)^8}.
\end{equation}
\end{theorem}
Note that for $p=1/2$ (\ref{testBIU}) yields
\begin{equation*}
\Delta(1/2) = \kappa(1/2) - \tau^2(1/2) = \frac {9 \times 16}{625} \left[ 4 + 16 - 20 + \frac {15}2 - \frac {15}{16} \right] = \frac{189}{125},
\end{equation*}
as it should.
Computation also shows that for $p=1/2$ (\ref{statBIU}) reduces to (\ref{Aemej}).

\begin{figure}[h]
\begin{center}
\includegraphics[height=3.5cm,width=10cm]{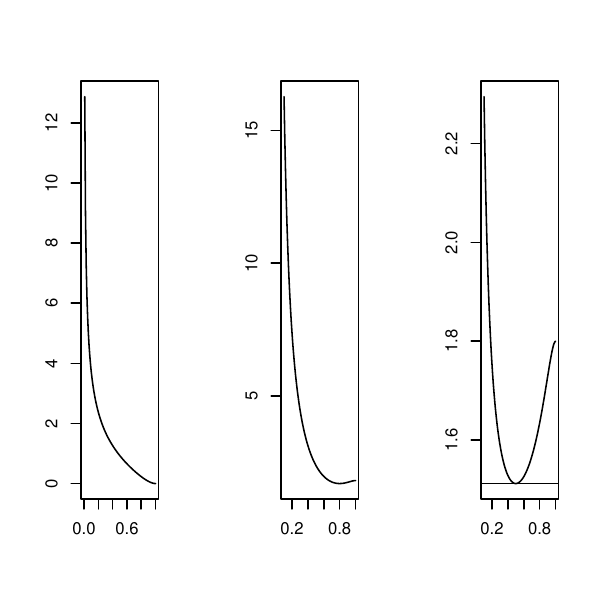}
\vspace{-0.8cm}
\caption{BIU($p$): Skewness, kurtosis and test function as functions of $p$.\label{figuurA}}
\end{center}
\end{figure}

Figure \ref{figuurA} gives plots of the skewness, kurtosis and the test function, i.e., kurtosis minus squared skewness, as functions of $p$.
Table \ref{tabelA} gives simulated power values ${\hat \beta}(p)$ based on 10,000 samples of size $n$ of the indicated BIU($p$) distribution.
The values lower and upper give the limits of the $95\%$ confidence interval around the estimated power ${\hat \beta}(p)$.

\begin{table}[h]
\begin{center}
\caption{BIU($p$): Simulated power values of the test (\ref{test}).\label{tabelA}}
\begin{tabular}{|c|c|c|c|c|c|c|c|c|c|}
  \hline
   & \multicolumn{3}{|c|}{$n$=50} &\multicolumn{3}{|c|}{$n$=100}&\multicolumn{3}{|c|}{$n$=1000}\\
  \hline
  $p$  & lower  & ${\hat \beta}(p)$  & upper & lower  & ${\hat \beta}(p)$  & upper & lower  & ${\hat \beta}(p)$  & upper     \\
  \hline
  0.4&  0.052&\bf 0.056& 0.061&  0.046 &\bf 0.051& 0.055&  0.032& \bf 0.036 &0.040 \\
  0.5&  0.047 & \bf 0.052& 0.056& 0.042& \bf 0.046& 0.051&  0.042 &\bf 0.046& 0.050\\
  0.6& 0.029& \bf 0.032& 0.036&0.019& \bf 0.022& 0.025&  0.010&\bf 0.012& 0.014 \\
  \hline
\end{tabular}
\end{center}
\end{table}

\section{Multimodality\label{multimodality}}

Let us consider the family of mixed normal distributions, denoted by MN($p,d$), that are mixtures of two normal distributions with equal variances 1 and means 0 and $d$, respectively.
The densities are given by
\begin{equation*}
f_{p,d}(x) = (1-p) \varphi(x) + p \varphi(x-d),\quad -\infty < x < \infty,
\end{equation*}
where $\varphi$ denotes the standard normal density.
To the results of \cite{Eisenberger} we add
\begin{lemma}\label{bimodality}
$f_{p,d}$ is bimodal if $(p \wedge(1-p)) \exp(d^2/8) \geq 1$ holds.
\end{lemma}
\noindent
{\bf Proof} \\
Note that the condition implies $f_{p,d}(0) > f_{p,d}(d/2) < f_{p,d}(d)$, which means bimodality.
\hfill
$\Box$ \\

\begin{figure}[ht]
\begin{center}
\includegraphics[height=3.5cm,width=4cm]{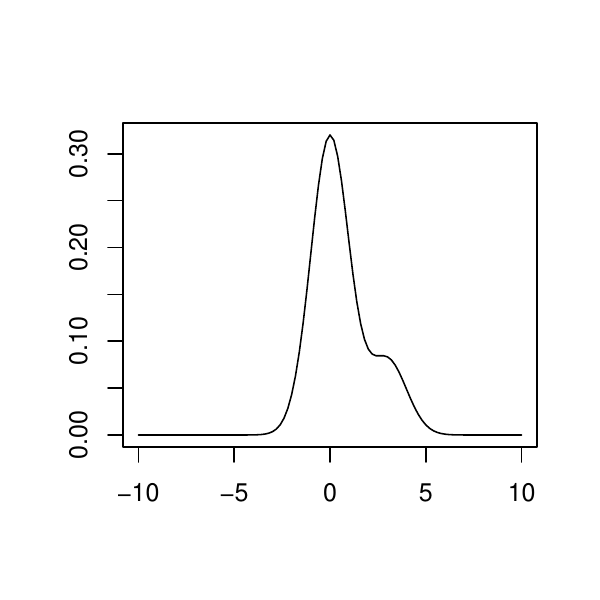}
\vspace{-0.8cm}
\includegraphics[height=3.5cm,width=4cm]{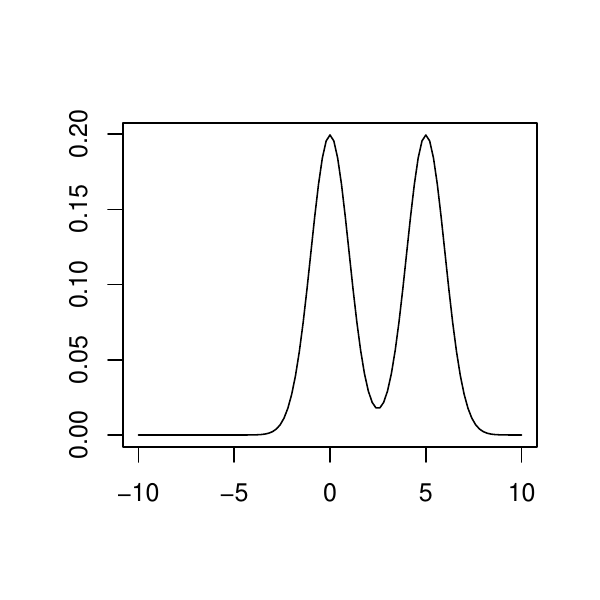}
\caption{Left: $p=0.2, d=3$. Right: $p=0.5, d=5$.\label{figuurC}}
\end{center}
\end{figure}

For $d=3$, the density changes from unimodal to bimodal at $p=0.2$, approximately.
For $p=0.5$ and $d=5$  the test function below is approximately equal to the bound $189/125$ (see Figure \ref{figuurC}) .

The following theorem gives the values for skewness, kurtosis and the test function for this distribution.
\begin{theorem}\label{MN}
For the {\rm MN($p,d$)} distribution the skewness is given by
\begin{equation*}
\tau(p,d)=d^3\,\frac{p(1-p)(1-2p)}{(1+d^2p(1-p))^{3/2}}.
\end{equation*}
The kurtosis is given by
\begin{equation*}
\kappa(p,d)=3+d^4\,\frac{p(1-p)(1-6p(1-p))}{(1+d^2p(1-p))^2}.
\end{equation*}
Finally the test function kurtosis minus squared skewness equals
\begin{eqnarray}\label{M1}
\lefteqn{ \Delta(p,d) = \kappa(p,d) - \tau^2(p,d) } \\
&& = 3 + d^4 p(1-p) \frac { 1 - 6 p(1-p) -2 d^2 p^2 (1-p)^2 } { (1+d^2 p(1-p))^3 }. \nonumber
\end{eqnarray}
\end{theorem}

\begin{figure}[h]
\begin{center}
\includegraphics[height=3.5cm,width=10cm]{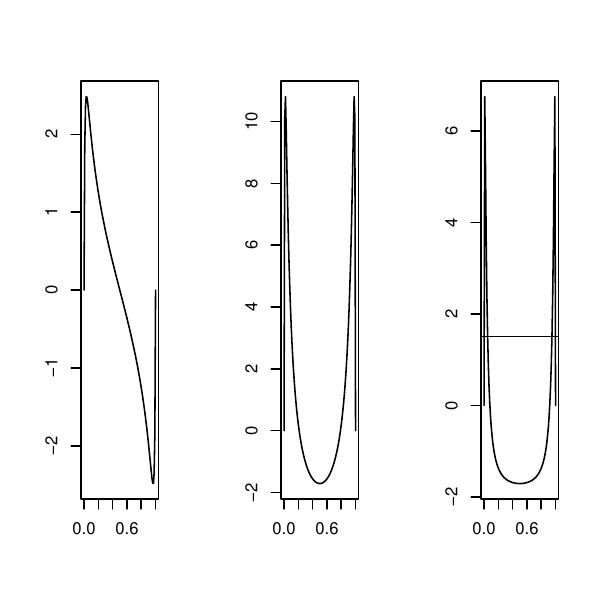}
\caption{MN($p,d)$: Skewness, kurtosis and test function as function of $p$ with $d=7$.\label{d=7}}
\end{center}
\end{figure}
\begin{table}[h]
\begin{center}
\caption{MN($p,d)$: Mixed normal, MN(0.5,$d$). Simulated  power values of the test (\ref{test}). \label{tabelB}}
\begin{tabular}{|c|c|c|c|c|c|c|c|c|c|}
  \hline
   & \multicolumn{3}{|c|}{$n$=50} &\multicolumn{3}{|c|}{$n$=100}&\multicolumn{3}{|c|}{$n$=1000}\\
  \hline
  $d$  & lower  & ${\hat \beta}(d)$  & upper & lower  & ${\hat \beta}(d)$  & upper & lower  & ${\hat \beta}(d)$  & upper     \\
  \hline
  5&  0.024 &\bf 0.027 &0.031&0.020 &\bf 0.023 &0.026& 0.015 &\bf 0.017& 0.020\\
 5.1 & 0.035&\bf  0.039& 0.043&0.036&\bf 0.040 &0.043& 0.061& \bf 0.066& 0.071\\
 5.2& 0.048&\bf  0.052& 0.057& 0.053&\bf 0.058& 0.062&  0.18&\bf 0.19 &0.20\\
  \hline
\end{tabular}
\end{center}
\end{table}

Table \ref{tabelB} gives simulated power values based on 10,000 samples of size $n$ of the indicated MN($p,d$) distribution.
The values lower and upper give the limits of the $95\%$ confidence interval around the estimated power ${\hat \beta}(d)$.

\begin{figure}[h]
\begin{center}
\includegraphics[height=9cm,width=9cm]{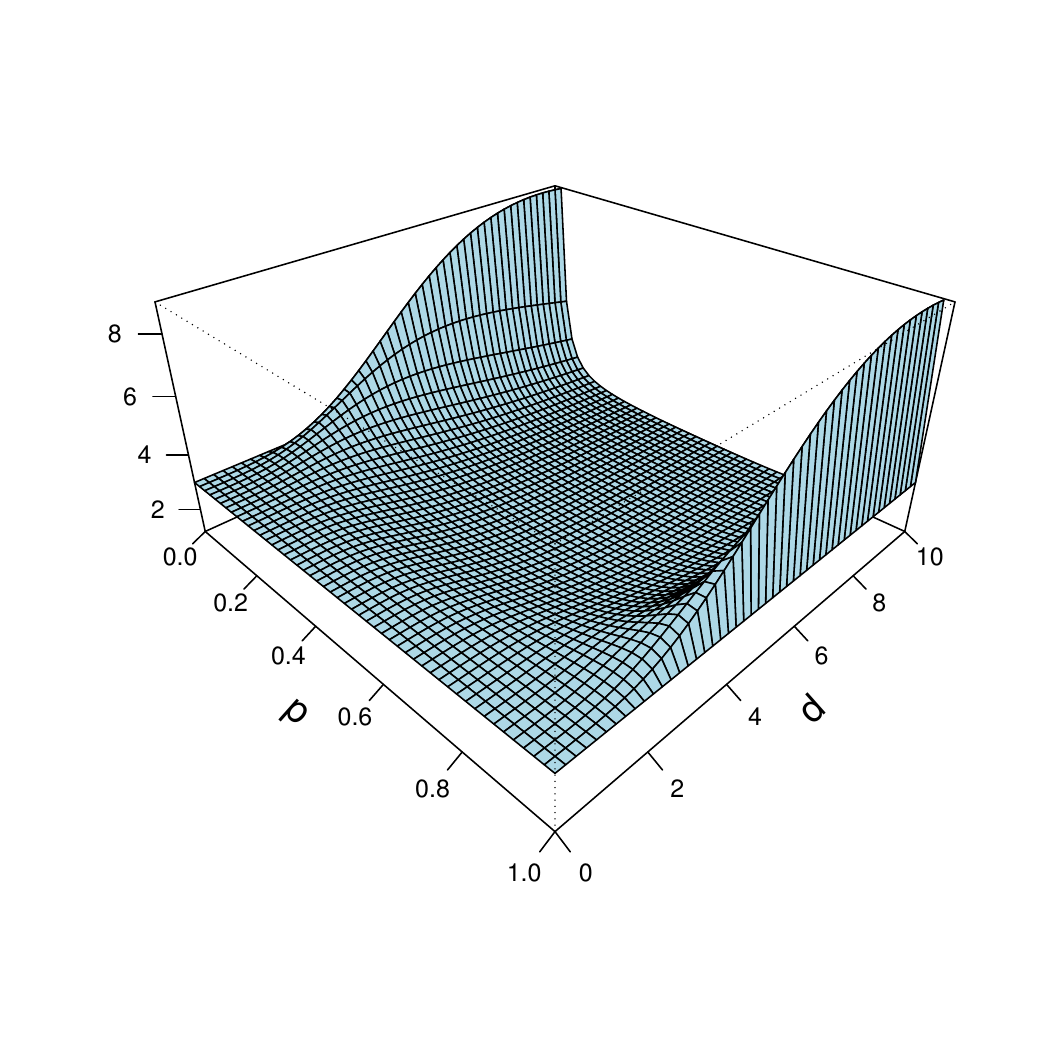}
\vspace{-1.5cm}
\caption{MN($p,d)$: A two dimensional plot of the test function for $p$ between zero and one and $d$ between zero and ten.\label{tweedimplot}}
\end{center}
\end{figure}

With $d = p^{-3/4}$ and $p$ small the test function equals approximately
\begin{equation*}
\Delta(p,p^{-3/4}) \approx 3 + \frac { p^{-2} } { p^{-3/2} } \left[ 1 - 6p - 2 p^{1/2} \right] \approx p^{-1/2},
\end{equation*}
which converges to $\infty$ as $p \downarrow 0$,
as is illustrated in Figure \ref{tweedimplot}.
Since also $p \exp(d^2/8)\\ > p,\ d^2/8 = p^{-1/2}/8$ holds, Lemma \ref{bimodality} shows that for $d = p^{-3/4}$ and $p$ small the mixed normal distribution is bimodal.
This means that for multimodal distributions the test function can be arbitrarily large.

Next we consider the following mixtures of two exponential densities
\begin{equation*}
g_{p,d}(x) = (1-p)e^{-x} + p e^{-x+d} {\bf 1}_{[x \geq d]}, \quad 0 < p < 1,\ 0 < d,\ 0 < x,
\end{equation*}
which are obviously bimodal.
If $X$ has density $g_{p,d}$, then computations shows
\begin{eqnarray*}
&& EX^k = (1-p)k! + p \sum_{j=0}^k {k \choose j} d^j (k-j)!
            = k! \left[ 1-p + p \sum_{j=0}^k \frac {d^j}{j!} \right], \\
&& EX = 1 + pd, \quad {\rm var~}X = 1 +d^2p(1-p), \\
&& E(X-EX)^3 = 2 + d^3p(1-p)(1-2p) , \\
&& E(X-EX)^4 = 9 + 6d^2p(1-p) +d^4p(1-p)(1-3p+3p^2), \\
&& \Delta(p,d) = \kappa(p,d) - \tau^2(p,d) \\
&& \quad \quad 5 - d^3p(1-p) \frac {4 - 8p - d (1-12p(1-p)) + 4 d^3 p^2 (1-p)^2 } { (1+d^2p(1-p))^3 }.
\end{eqnarray*}
With $d=p^{-1/3}$ we see that for $p \downarrow 0$ we have
\begin{equation}\label{Deltalarge}
\Delta(p,p^{-1/3}) \approx 5  -4 + p^{-1/3} \approx p^{-1/3}.
\end{equation}
This shows that for the choice of the distance $d=p^{-1/3}$ between the modes the test function can be made arbitrarily large by choosing $p$ sufficiently small.

\section{Skewness-kurtosis sets}\label{sksets}

In this section we present the exact shape of three skewness-kurtosis sets.

\subsection{The skewness-kurtosis set for arbitrary distributions}\label{sk}

\begin{theorem}\label{sks1}
For every point $(x_0,y_0) \in \{ (x,y) \in {\mathbb R}^2 \, | \, y \geq x^2 +1 \}$ there exists a distribution whose skewness equals $x_0$ and whose kurtosis equals $y_0$.
\end{theorem}
\noindent
{\bf Proof} \\
Consider the following well defined random variable
\begin{equation}
X = \left\{  \begin{array}{lcl}
            - \sqrt{\alpha(1+z)}       &                          & \frac 1{(1+\alpha)(1+z)} \\
            0                          & {\rm ~with~probability~} & 1 - \frac 1{1+z}         \\
            \sqrt{ \frac {1+z}\alpha } &                          & \frac \alpha {(1+\alpha)(1+z)}
            \end{array}  \right.
\end{equation}
with $0 < \alpha \leq 1$ and $0 \leq z$.
Some computation shows
\begin{eqnarray*}
\lefteqn{ EX = 0, \quad EX^2 = 1, \quad \tau = EX^3 = \frac {1-\alpha}{\sqrt \alpha} \sqrt{1+z}, } \\
 &&  \hspace{5em} \kappa = EX^4 = \frac {1-\alpha + \alpha^2}\alpha (1+z). \hspace{5em}.
\end{eqnarray*}
Note that the function $\psi(\alpha) = (1-\alpha)/\sqrt{\alpha},\ 0 < \alpha \leq 1,$ is continuous with $\psi(1) = 0$ and
$\lim_{\alpha \downarrow 0} \psi(\alpha)= \infty$.
This means that $\tau$ can take all values in $[0,\infty)$.
Furthermore, we have $\kappa = \tau^2 +1 +z$ and hence for appropriate choices of $\alpha$ and $z$ each point in the set
$\{ (x,y) \in {\mathbb R}^2 \, | \, 0 \leq x,\ y \geq x^2 +1 \}$ can be attained by the corresponding $(\tau, \kappa)$.
Finally, with $-X$ the set $\{ (x,y) \in {\mathbb R}^2 \, | \, x \leq 0,\ y \geq x^2 +1 \}$ can be attained.
We have shown that for every point on or above the parabola $y = x^2 +1$ there exists a distribution such that its skewness and kurtosis coincide with this point.
\hfill
$\Box$ \\

\subsection{The skewness-kurtosis set for unimodal distributions}\label{skunimodal}

\begin{theorem}\label{sks2}
Define
\begin{eqnarray*}
\lefteqn{ L(\alpha , \tau) = \frac 15 \left(4+\alpha^2 \right)^{-3}
\left\{ \frac 13 \left(1+4\alpha^2 \right) \left(4+\alpha^2 \right)^3 \tau^2 \right. } \\
&& \quad + 2 \alpha \left( 4 + \alpha^2 \right)^{3/2} \left[ 8\sqrt{3} \alpha^4 + (2 \sqrt{3} +1)\alpha^2 -1\right] \tau \\
&& \quad \left. + 3 \left[48 \alpha^8 +(15 + 4 \sqrt{3} ) \alpha^6 + (60-4 \sqrt{3}) \alpha^4 +60 \alpha^2 + 192 \right] \right\}.
\end{eqnarray*}
The kurtosis $\kappa$ of a unimodal distribution with skewness $\tau$ can take any value at least equal to
\begin{equation}\label{sks2.1}
\tau^2 + \min_{\alpha \in {\mathbb R}} L(\alpha , \tau).
\end{equation}
\end{theorem}
\begin{figure}[h]
\begin{center}
\includegraphics[height=9cm,width=9cm]{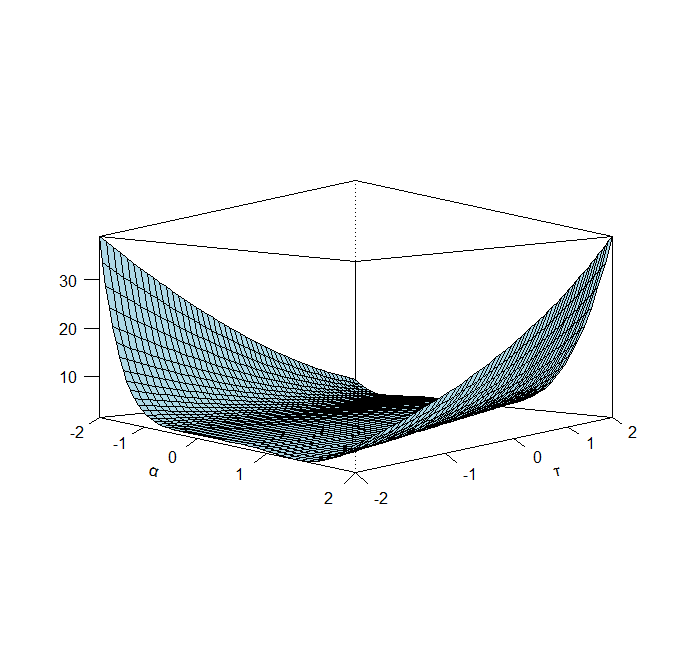}
\vspace{-2cm}
\caption{A two dimensional plot of the function $L$.\label{Lplot1}}
\end{center}
\end{figure}

\noindent
{\bf Proof} \\
By Khintchine's representation a unimodal random variable $X$ with mode at 0 may be written as
$X = UZ$ with $U$ uniform on $(0,1)$ and $U$ and $Z$ independent; cf. Theorem V.9, p. 158, of Feller (1971) or Lemma A.1 of \cite{Ion}.
With $G$ the distribution of $Z$ we write $\mu_G, \sigma_G, \tau_G$ and $\kappa_G$ for mean, standard deviation, skewness,
and kurtosis, respectively, of $G$.
Tedious computation shows that the skewness $\tau$ and kurtosis $\kappa$ of $X$ satisfy (cf. (11) of \cite{Mokveld})
\begin{eqnarray*}
\lefteqn{ \kappa = \tau^2 -108 \left(4 + \left( \frac{\mu_G}{\sigma_G} \right)^2 \right)^{-3} \left( \tau_G + \frac{\mu_G}{\sigma_G} \right)^2 } \\
&& + \frac 95 \left(4 + \left( \frac{\mu_G}{\sigma_G} \right)^2 \right)^{-2}
\left( 16 \kappa_G + 24 \frac{\mu_G}{\sigma_G}\, \tau_G +16 \left( \frac{\mu_G}{\sigma_G} \right)^2 + \left( \frac{\mu_G}{\sigma_G} \right)^4 \right).
\end{eqnarray*}
For every $z \geq 0$ there exists a distribution $G$ with $\kappa_G = \tau_G^2 +1 + z$  in view of Theorem \ref{sks1}.
Since skewness and kurtosis are location and scale invariant we may assume $\mu_G / \sigma_G = \alpha$ for any $\alpha \in \mathbb{R}$.
Computation shows
\begin{equation*}
\tau_G = \alpha^3 + \frac 1{6 \sqrt{3}} \left( 4 + \alpha^2 \right)^{3/2} \tau.
\end{equation*}
With these substitutions we arrive at
\begin{equation*}
\kappa = \tau^2 + L(\alpha , \tau) + \frac{144}5 \left( 4 + \alpha^2 \right)^{-2} z,
\end{equation*}
which implies the theorem.
\hfill
$\Box$ \\

According to Theorem 1.1 of \cite{Mokveld} equality holds in (\ref{Philip}) if and only if $X$ has a one-sided boundary-inflated uniform distribution with mass 1/2 at the atom.
Computation shows that such a distribution has skewness and kurtosis equal to
\begin{equation}\label{sks2.1}
\tau = \pm \frac 65 \, \sqrt{\frac 35}, \quad \kappa = \frac {297}{125}.
\end{equation}
This means that the parabola $\{ (x,y) \in {\mathbb R}^2 \, | \, y = x^2 +189/125 \}$ contains only two points of the skewness-kurtosis set of unimodal distributions.
With $|\alpha| = 1$ and $\tau$ as in (\ref{sks2.1}) such that $\alpha \tau$ is negative, straightforward computation gives
$L(\alpha , \tau) = 189/125$ in agreement with Theorem 1.1 of \cite{Mokveld}.

\begin{figure}[ht]
\begin{center}
\includegraphics[height=4cm,width=4cm]{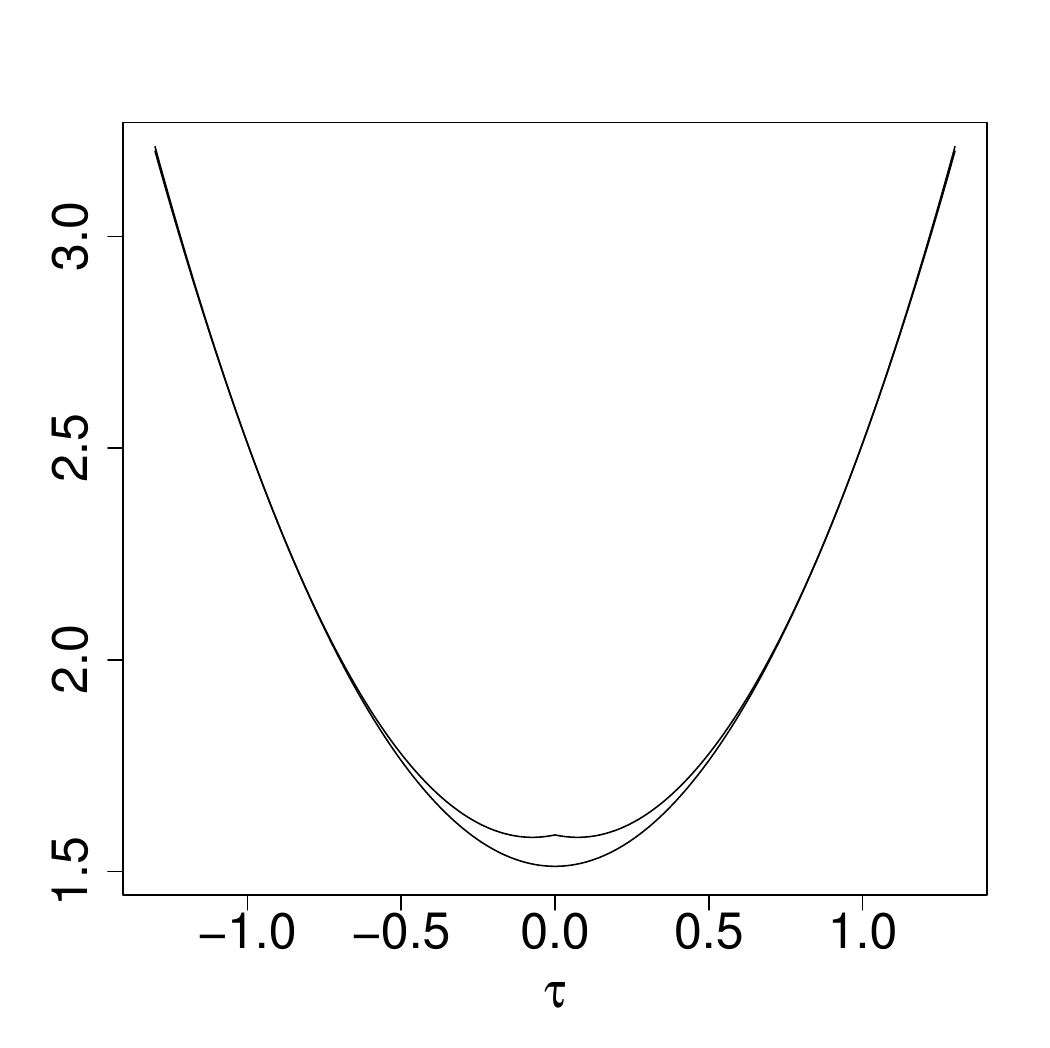}
\includegraphics[height=4cm,width=4cm]{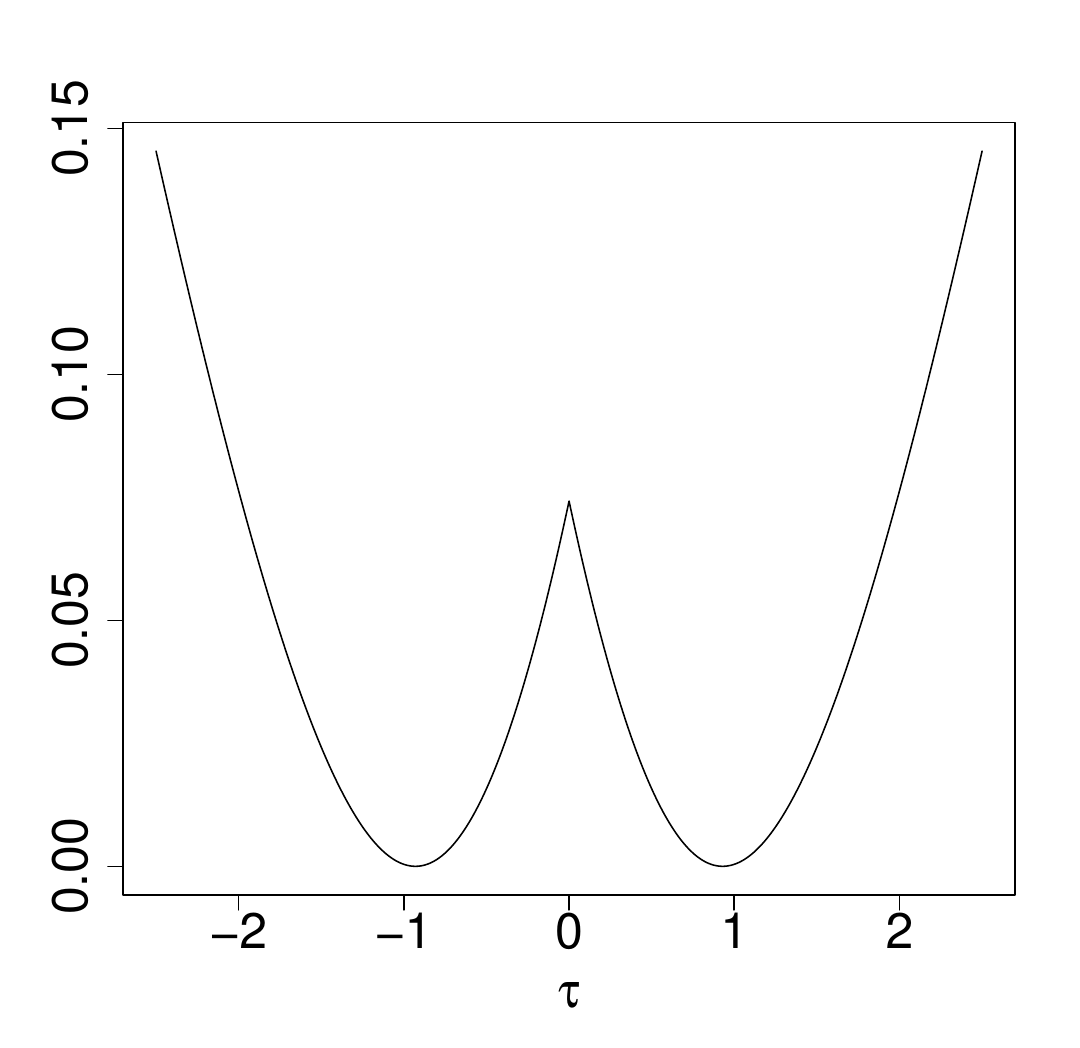}
\caption{Left: The two functions $f(\tau)$ and $b(\tau)$. Right: Their difference. \label{boundsplots}}
\end{center}
\end{figure}

Let us denote the boundary function given by (\ref{sks2.1}) by $b(\tau)$ and the boundary function for $\kappa$ in \cite{Mokveld} by $f(\tau)=\tau^2+189/125$ .
The function $b(\tau)$ can be computed by numerical minimization. For unimodal distributions we have $b(\tau)\geq f(\tau)$ with, as argued above, equality for $\tau=\pm 2 (3/5)^{3/2}\approx  \pm 0.9295$. This is illustrated by Figure \ref{boundsplots}.

\subsection{The skewness-kurtosis set for symmetric unimodal distributions}\label{sksunimodal}

\begin{theorem}\label{sks3}
For symmetric unimodal distributions the skewnes-kurtosis set equals $\{ (0,y) \in {\mathbb R}^2 \, | \, y \geq 9/5 \}$.
\end{theorem}
\noindent
{\bf Proof} \\
In view of (\ref{Bert}) it remains to be shown that the kurtosis of a symmetric unimodal distribution can take any value at least equal to $9/5$.
Let $X$ be a random variable that is uniformly distributed on the interval $[-1,1]$ with probability $p$ and has an atom at 0 with mass $1-p$.
The $k$-th moment of $X$ equals
\begin{equation*}
E(X^k) = \frac p{k+1} \ {\bf 1}_{[k {\rm ~even}]}, \quad k=1, 2, \dots,
\end{equation*}
and hence we have
\begin{equation*}
\mu = 0, \quad \sigma^2 = \frac p3, \quad \tau = 0, \quad \kappa = \frac 9{5p}.
\end{equation*}
Since the map $(0,1] \to [9/5, \infty), \ p \mapsto 9/(5p),$ is onto, the proof is complete.
\hfill
$\Box$ \\

\begin{figure}[ht]
\begin{center}
\includegraphics[height=4cm,width=4cm]{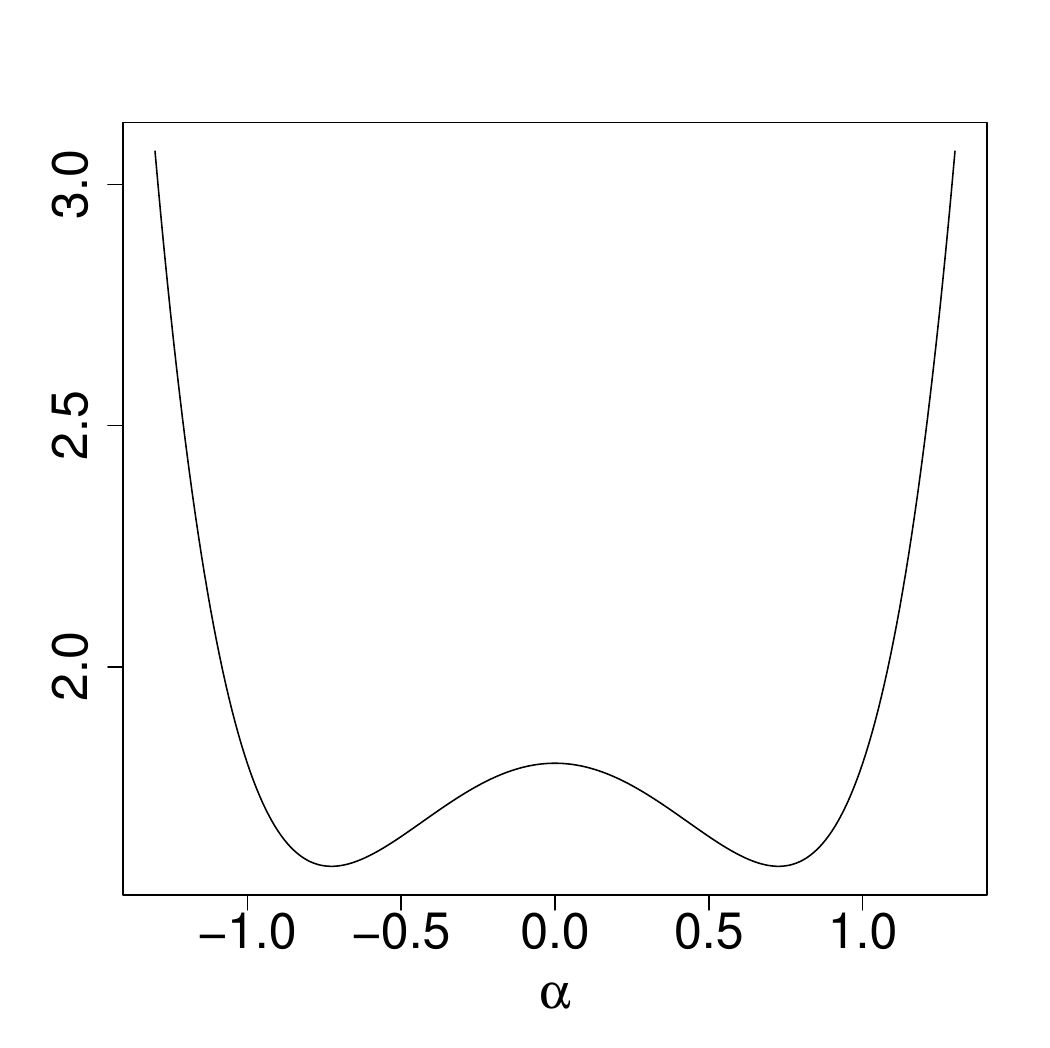}
\caption{The function $L(\alpha,0)$. \label{Ltau=0}}
\end{center}
\end{figure}

By the proof of Theorem \ref{sks2} we obtain for unimodal distributions with skewness $\tau = 0$
\begin{equation}\label{sks3.1}
\kappa \geq \min_{\alpha \in {\mathbb R}} L(\alpha , 0).
\end{equation}
Since computation shows $L(-1,0)=L(0,0)=L(1,0)=9/5$, Figure \ref{Ltau=0} implies that the lower bound in (\ref{sks3.1}) is less than $9/5$.
Together with Theorem \ref{sks2} we conclude that there exist nonsymmetric unimodal distributions with skewness $\tau=0$ and $\kappa < 9/5$.

\section{Proofs}\label{proofs}

\subsection{Proof of Theorem \ref{BIU}}

Let $X=BU$ be a random variable with $B$ and $U$ independent, $B$ Bernoulli distributed with success probability $p$ and $U$ uniform on the unit interval.
The moments and variance of $X$ are
\begin{eqnarray*}
&& EX^k = EB\, EU^k = \frac p{k+1},\ k \neq 0, \\
&& {\rm var~} X = \frac p3 - \left( \frac p2 \right)^2 = \frac p{12}(4-3p).
\end{eqnarray*}
This implies
\begin{eqnarray*}
\lefteqn{ \tau^2(p) = ({\rm var~}X)^{-3} \left( E(X - EX)^3 \right)^2 } \\
&& = \left( \frac {12}{p(4-3p)} \right)^3 \left( \frac p4 -3 \frac p3 \, \frac p2 + 2 \left( \frac p2 \right)^3 \right)^2 \\
&& = \frac {108}{p(4-3p)^3} (1-2p+p^2)^2 = \frac {108}{p(4-3p)^3} (1-p)^4
\end{eqnarray*}
and
\begin{eqnarray*}
\lefteqn{ \kappa(p) = ({\rm var~}X)^{-2}  E(X - EX)^4 } \\
&& = \left( \frac {12}{p(4-3p)} \right)^2 \left( \frac p5 -4 \frac p4 \, \frac p2 + 6 \frac p3  \left( \frac p2 \right)^2
     -3 \left( \frac p2 \right)^4 \right) \\
&& = \frac {144}{5p(4-3p)^2} \left(1 - \frac 52 p + \frac 52 p^2 - \frac{15}{16} p^3 \right) \\
&& = \frac {9}{5p(4-3p)^2} \left(16 - 40 p + 40 p^2 - 15 p^3 \right).
\end{eqnarray*}
Consequently, we have
\begin{equation*}
\Delta(p) = \kappa(p) - \tau^2(p) = \frac 9{5p(4-3p)^3} \left[ 4 + 32p - 80 p^2  + 60 p^3 - 15 p^4 \right].
\end{equation*}

Furthermore, we have for $k = 0, 1, \dots$
\begin{eqnarray*}
\lefteqn{ E(X-EX)^k = \sum_{j=0}^k {k \choose j} (-\tfrac 12 p)^j EX^{k-j}
= (-\tfrac 12 p)^k + \sum_{j=0}^{k-1} {k \choose j} (-\tfrac 12 p)^j \frac p{k-j+1} } \\
&& \hspace{4em} = (-\tfrac 12 p)^k + \frac p{k+1} \sum_{j=0}^{k-1} {{k+1} \choose j} (-\tfrac 12 p)^j  \\
&& \hspace{4em} = (-\tfrac 12 p)^k + \frac p{k+1} \left[ (1-\tfrac 12 p)^{k+1} - (k+1)(-\tfrac 12 p)^k - (-\tfrac 12 p)^{k+1} \right]. \nonumber
\end{eqnarray*}
Consequently the $k$th reduced moment equals
\begin{eqnarray}\label{b7}
\lefteqn{ \nu_k = E \left( \frac {X-EX} {\sqrt{{\rm var~}X}} \right)^k
= \left( \frac p{12}(4-3p) \right)^{-k/2} \nonumber } \\
&& \times \left( (-\tfrac 12 p)^k + \frac p{k+1} \left[ (1-\tfrac 12 p)^{k+1} - (k+1- \tfrac 12 p)(-\tfrac 12 p)^k \right] \right).
\end{eqnarray}
Straightforward but tedious computation yields (\ref{statBIU}).

\subsection{Proof of Theorem \ref{MN}}

Let $X=N+dB$ be a random variable with $B$ and $N$ independent, $B$ Bernoulli distributed with success probability $p$ and $N$ standard normal.
The moments of $X$ can be computed by
\begin{eqnarray*}
EX^k &=& E(N+dB)^k=E\sum_{j=0}^k{k \choose j}N^jd^{k-j}B^{k-j}\nonumber\\
&=&\sum_{j=0}^{k-1}d^{k-j}{k \choose j}EN^jEB+EN^k=EN^k+p\sum_{j=0}^{k-1}d^{k-j}{k \choose j}EN^j.
\end{eqnarray*}
This implies that the expectation and variance of $X$ are
\begin{equation*}
EX=dp,\quad EX^2=1+d^2 p,\quad {\rm var} X=1+d^2 p(1-p).
\end{equation*}
Furthermore, we have
\begin{eqnarray*}
\lefteqn{\tau^2(p,d) = ({\rm var~}X)^{-3} \left( E(X - EX)^3 \right)^2 \nonumber } \\
&& = ( 1+d^2p(1-p)))^{-3} \left(d^3(p-3p^2+2p^3) \right)^2 = d^6\,\frac{p^2(1-p)^2(1-2p)^2}{(1+d^2p(1-p))^3},
\end{eqnarray*}
and
\begin{eqnarray*}
\kappa(p,d) = ({\rm var~}X)^{-2}  E(X - EX)^4 =  3 + d^4 p(1-p) \frac { 1-6p(1-p) } { (1 + d^2 p(1-p) )^2 },
\end{eqnarray*}
which combined yield (\ref{M1}).

\end{document}